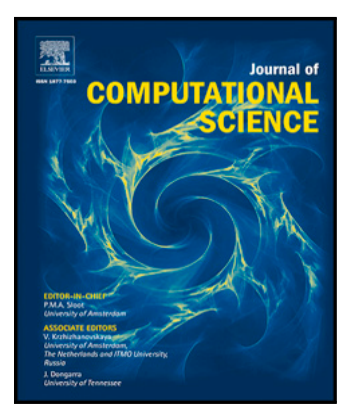

# Some observations regarding the RBF-FD approximation accuracy dependence on stencil size

Andrej Kolar-Požun [a,b,*], Mitja Jančič [a], Miha Rot [a,c], Gregor Kosec [a]

[a] Jožef Stefan Institute, Parallel and Distributed Systems Laboratory, Jamova 39, 1000 Ljubljana, Slovenia
[b] University of Ljubljana, Faculty of Mathematics and Physics, Jadranska 19, 1000 Ljubljana, Slovenia
[c] International Postgraduate School Jožef Stefan, Jamova 39, 1000 Ljubljana, Slovenia

ABSTRACT

When solving partial differential equations on scattered nodes using the Radial Basis Function-generated Finite Difference (RBF-FD) method, one of the parameters that must be chosen is the stencil size. Focusing on Polyharmonic Spline RBFs with monomial augmentation, we observe that it affects the approximation accuracy in a particularly interesting way — the solution error oscillates under increasing stencil size. We find that we can connect this behaviour with the spatial dependence of the signed approximation error. Based on this observation we are able to introduce a numerical quantity that could indicate whether a given stencil size is locally optimal. This work is an extension of our ICCS 2023 conference paper (Kolar-Požun et al., 2023).

## 1. Introduction

Partial Differential Equations (PDEs) are ubiquitous in science and engineering, as well as economics and related fields. For this reason, study of PDEs and their solutions is an extremely active area of research. As many real-world cases are too complicated to solve analytically, a substantial part of this research focuses on numerical PDE solution procedures.

Currently, the most widely used method of obtaining a numerical solution to a given PDE is the finite element method (FEM) [1]. It is well understood, supports all types of adaptivities [2] accompanied with well understood error indicators [3] and isogeometric analysis [4]. A necessary step in FEM analysis is meshing, where the entire computational domain is covered with polyhedrons. In realistic 3D domains the process of meshing often requires user assistance and thus cannot be automated [5]. This motivates research into meshless methods [6] that work directly on point clouds often referred to as "scattered nodes", which considerably simplifies the discretisation of the domain regardless of its dimensionality or shape [7,8]. The elegant formulation of meshless methods is also convenient for implementing h- [9] and hp- [10] adaptivities, considering different approximations of partial differential operators in terms of the stencil shape and size [11,12] and the local approximation order [13].

Over the years numerous meshless methods were proposed, some of the commonly used ones being the Finite Point Method [14], Smooth Particle Hydrodynamics [15], the Generalised Finite Difference Method [16] and the Moving Least Squares [17] to name a few. In this paper we focus on a yet another popular choice - a method based on Radial Basis Functions (RBFs). RBFs have first seen widespread use in interpolation on unstructured discretisations in the 1970s [18] and have become attractive due to their positive definiteness property [19], resulting in a well-defined interpolation problem for any pairwise distinct collection of nodes. In 1990, Kansa used the RBF interpolants to also approximate linear differential operators, resulting in a RBF-based collocation method for solving PDEs [20].

In the early 2000s, a local variant of Kansa's method was proposed, with improved stability properties and lower computational complexity [21,22]. This method, today known as the Radial Basis Function-generated Finite Difference (RBF-FD) method, approximates a linear differential operator only locally, given a chosen neighbourhood. The neighbourhood in question is referred to as a stencil of a given point and is commonly chosen to simply consist of its $n$ closest points.

In the last few years, several improvements of the RBF-FD method have been proposed. For instance, RBF-HFD (Hermite Finite Difference) accounting for derivative information for added accuracy [23], oversampled RBF-FD with increased stability properties [24] and overlapped RBF-FD with improved computational complexity [25]. However, for the purposes of this paper we limit our analyses to the simplest and most popular original version [21].

Initially, the researchers in this field focused on smooth RBFs such as Gaussians or Inverse Multiquadrics [26,27]. Their downside is that they all possess a free parameter — the shape parameter, which can substantially affect the method's accuracy and the precise interpretation of which is still a major subject of research, making it difficult to choose appropriately [22].

* Corresponding author at: Jožef Stefan Institute, Parallel and Distributed Systems Laboratory, Jamova 39, 1000 Ljubljana, Slovenia.
*E-mail addresses:* andrej.pozun@ijs.si (A. Kolar-Požun), mtja.jancic@ijs.si (M. Jančič), miha.rot@ijs.si (M. Rot), gregor.kosec@ijs.si (G. Kosec).

For this reason, RBFs without a shape parameter have been gaining traction lately. Concretely, Polyharmonic Splines (PHS) with monomial augmentation have become a popular choice, becoming a subject of extensive research and applications [10,28–30].

The main topic of this paper is a study of the effect of stencil size $n$ on the quality of differential operator approximation in PHS RBF-FD.

Regarding stability, Bayona already noticed that increasing the stencil size smoothens the cardinal functions and can aid with the diagonal dominance of the resulting differentiation matrix [31].

As for approximation accuracy, for the case of smooth RBFs, it was observed that the stencil size can affect the order of the approximation [32]. On the other hand, in the context of PHS RBF-FD, it is easy to see that the order of the method is determined by the degree of augmented monomials [33]. The effect of the stencil size on approximation accuracy, however, is not so evident. The fact that stencil properties can affect the accuracy is sometimes mentioned in passing [24] but, to the best of our knowledge, has not been analysed in detail.

We observe that the choice of an appropriate stencil size in PHS RBF-FD can have a substantial impact on the accuracy. Moreover the accuracy of the method displays an oscillatory behaviour under increasing stencil size. In the remainder of the paper, we present our findings in greater detail. Our main objectives are determining why and under which conditions the oscillatory behaviour of the solution error occurs and whether we can predict its quantitative properties a priori. Namely, it would be beneficial if we could predict the optimal stencil size for a given problem and thus improve the method's accuracy without changing our discretisation set or increasing the order of our method.

The following section describes our problem setup along with the numerical solution procedure. In Section 3 our results are then discussed, where we start by demonstrating the aforementioned error oscillations. We then proceed by showing how these oscillations can be connected to the spatial dependence of pointwise error. In Section 4 several further analyses are performed. First, we check to what extent the oscillations remain if some aspect of the method is modified and then analyse how our observations change when considering different problem setups. In Section 5 a realistic problem of determining the stationary temperature profile of a heatsink is analysed as a potential application of our findings.

This paper is an extension of our conference paper [34], in which the aforementioned error oscillations were first observed.

## 2. Problem setup

Our analyses are performed on the case of the Poisson equation

$$\nabla^2 u(\mathbf{x}) = f(\mathbf{x}), \tag{1}$$

where the domain is an open disc $\Omega = \{\mathbf{x} \in \mathbb{R}^2 : \|\mathbf{x} - (0.5, 0.5)\| < 0.5\}$. We choose the function $f(\mathbf{x})$ such, that the problem given by Eq. (1) has a known solution. Concretely, we choose

$$u(x, y) = \sin(\pi x)\sin(\pi y), \tag{2}$$

$$f(x, y) = -2\pi^2 \sin(\pi x)\sin(\pi y) \tag{3}$$

with the Dirichlet boundary conditions given by a restriction of $u(\mathbf{x})$ to the boundary $\partial\Omega$.

We discretise the domain with the internodal distance $h = 0.01$, first discretising the boundary and then the interior using an advancing front method-based algorithm proposed in [8,35] that guarantees minimal spacing and local regularity of the resulting discretisation set. From here on out we will refer to it as the DIVG (Dimension Independent Variable density node Generation) algorithm. On the resulting discretisation points we then obtain a numerical solution to the Poisson problem using the RBF-FD procedure, which we will now briefly describe. For a more complete description we refer the reader to [36]. We associate to each discretisation point $\mathbf{x}_i$ its stencil, which consists of its $n$ nearest neighbours. Fig. 1 serves as an visual representation of how a typical stencil looks like. The Laplacian will be discretised locally on each such stencil. For readability, we will enumerate the nodes in a given stencil as $\mathbf{x}'_1, \dots \mathbf{x}'_n$. Let us start by assuming we want to interpolate an unknown function $u(\mathbf{x})$ given its values in the stencil nodes. We write an RBF interpolant with monomial augmentation as:

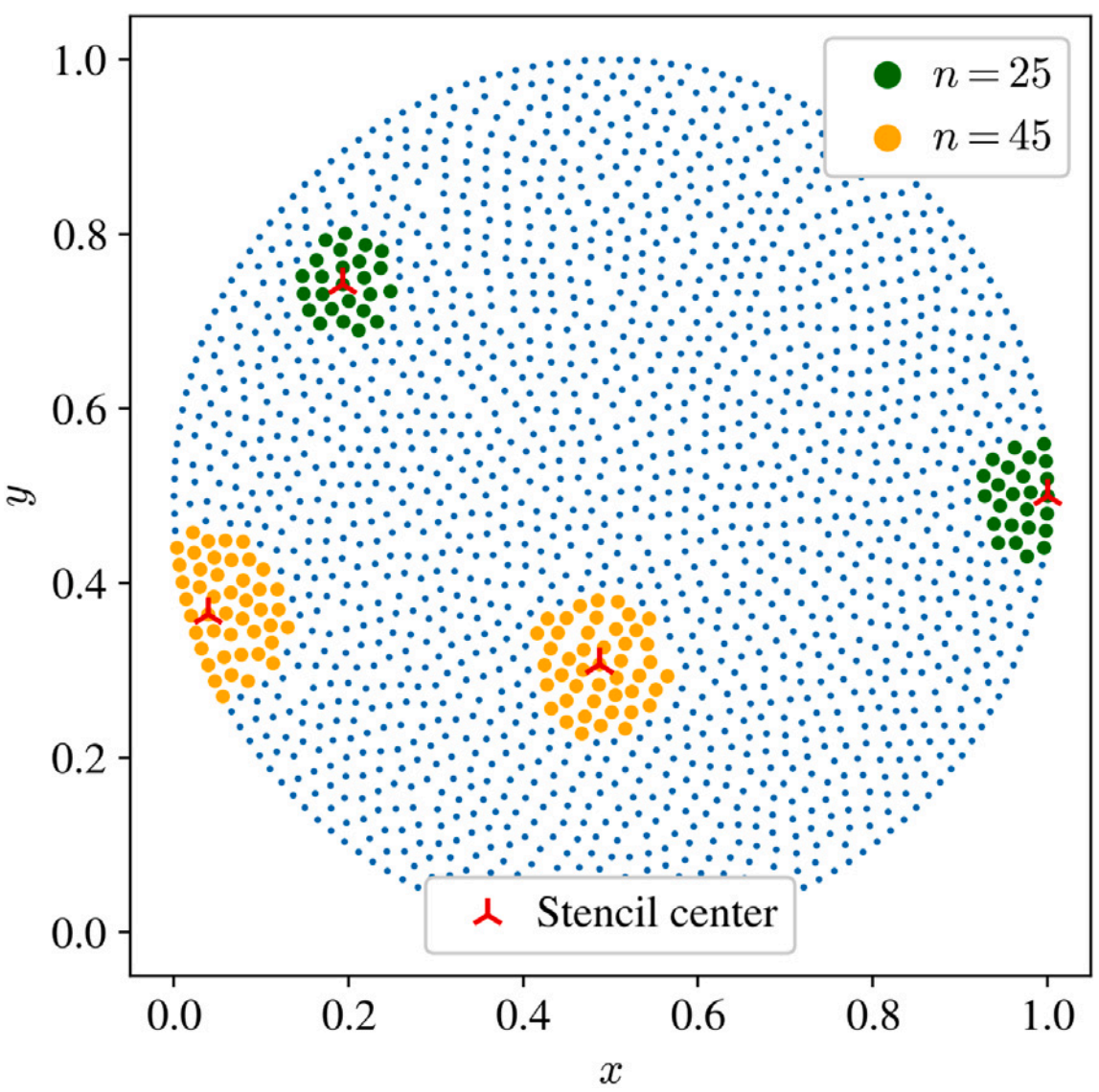

**Fig. 1.** Example discretisation set generated by the DIVG algorithm. Example stencils are also displayed.

$$s(\mathbf{x}) = \sum_{i=1}^{n} \alpha_i \phi(\|\mathbf{x} - \mathbf{x}'_i\|) + \sum_{i=1}^{M_m} \beta_i p_i(\mathbf{x}), \tag{4}$$

where the first sum is over the chosen RBFs. In our analyses we have opted for radial cubics $\phi(r) = r^3$. Second sum ranges over the $M_m$ basis functions of polynomials in two variables up to degree $m$, inclusive — concretely, this implies $M_m = (m+1)(m+2)/2$. In practice (and also in this paper) $p_i(\mathbf{x})$ are taken to be monomials. Unless otherwise stated, we will work with $m = 3$. Before proceeding, it will be useful to rewrite the interpolant in cardinal (also known as Lagrange) form:

$$s(\mathbf{x}) = \sum_{i=1}^{n} \psi_i(\mathbf{x}) u(\mathbf{x}'_i), \tag{5}$$

where $\psi_i(\mathbf{x})$ are known as the cardinal functions. Applying a chosen operator $\mathcal{L}$ (in our case $\nabla^2$), we can express:

$$\mathcal{L}(s(\mathbf{x})) = \sum_{i=1}^{n} \mathcal{L}(\psi_i(\mathbf{x})) u(\mathbf{x}'_i). \tag{6}$$

We will be interested in the operator value at the centre point of the stencil $x_C$:

$$\mathcal{L}(s(\mathbf{x}_C)) = \sum_{i=1}^{n} w_i u(\mathbf{x}'_i), \tag{7}$$

where we have defined $w_i = \mathcal{L}(\psi_i(\mathbf{x}_C))$. This gives us the differentiation weights approximating the operator $\mathcal{L}$ and directly generalising the finite difference method [30]. Skipping some details, more thoroughly explained in [36], let us just state that in practice the weights $w_i$ are obtained as a solution of an appropriate linear system. Having calculated the differentiation weights, we can now convert the PDE into a global (and sparse) linear system, similarly to any traditional finite difference methods. We then solve it to obtain an approximate solution $\hat{u}(\mathbf{x})$. The source code for the above solution procedure and also all the forthcoming analyses is readily available in our git repository[1]. Note that the provided code relies heavily on an open source meshless library Medusa [37].

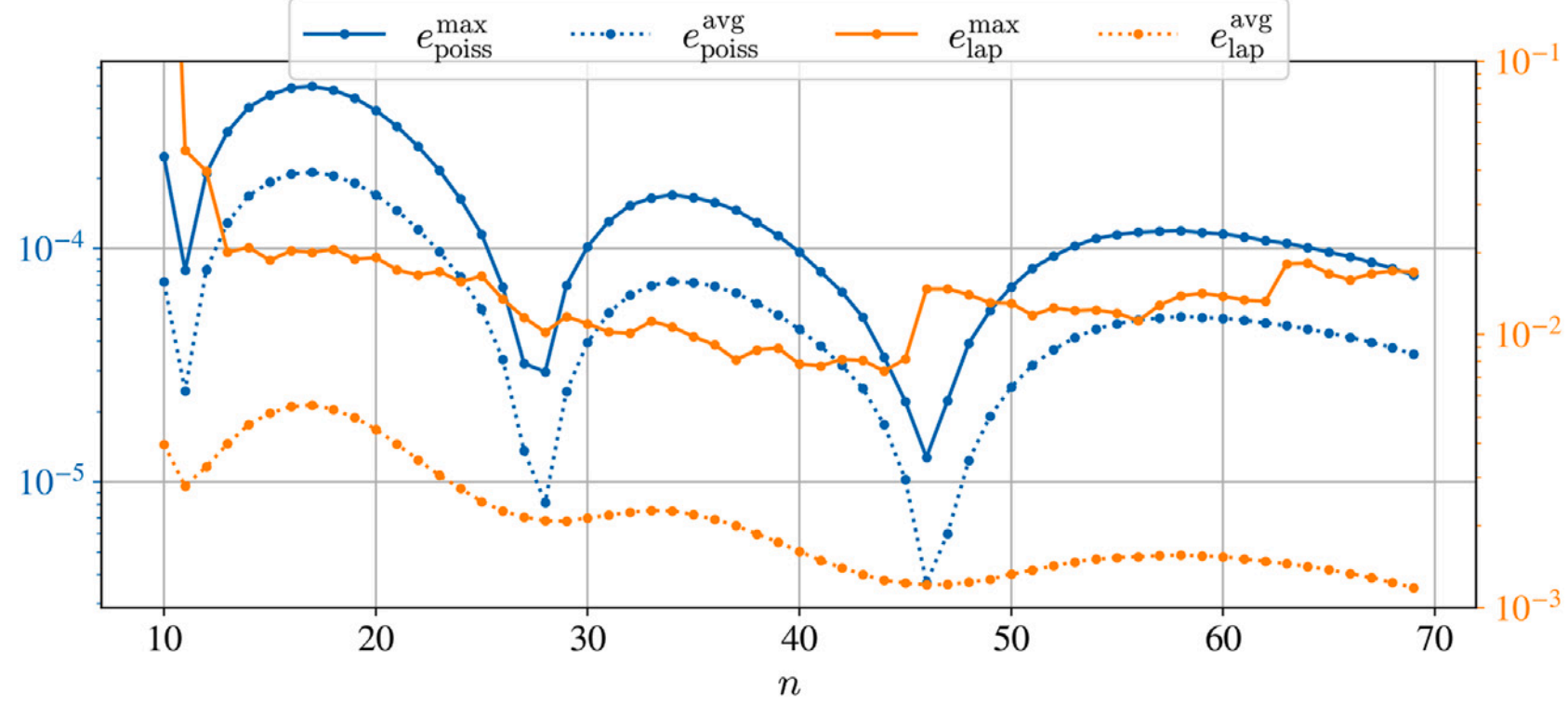

**Fig. 2.** Dependence of the approximation error on the stencil size $n$. The data used to make this plot is available in the appendix.

Regarding the described RBF-FD method, an important property is that the degree of monomial augmentation $m$ determines the approximation order — truncation error scales as $h^{m+1-l}$, where $l$ is the order of the differential operator $\mathcal{L}$ (in our case, $l = 2$) [24]. Compared with some of the simpler methods of same order, such as a purely monomial approximation, RBF-FD has improved stability properties. In fact, we can increase the method stability by keeping the monomial augmentation degree constant and increasing the stencil size [30].

Having both the analytical and approximate solutions, we will be interested in the approximation error. It will turn out to be useful to consider the signed pointwise errors of both the solution and the Laplacian approximation:

$$e^{\pm}_{\text{poiss}}(\mathbf{x}_i) = \hat{u}_i - u_i, \tag{8}$$

$$e^{\pm}_{\text{lap}}(\mathbf{x}_i) = \tilde{\nabla}^2 u_i - f_i, \tag{9}$$

where $\tilde{\nabla}^2$ is the discrete approximation of the Laplacian and we have introduced the notation $u_i = u(\mathbf{x}_i)$. The “poiss” and “lap” subscripts may be omitted in the text when referring to both errors at once.

As a quantitative measure of the approximation quality, we will also look at the average/max absolute value error:

$$e^{\max}_{\text{poiss}} = \max_{\mathbf{x}_i \in \Omega} |e^{\pm}_{\text{poiss}}(\mathbf{x}_i)|, \tag{10}$$

$$e^{\text{avg}}_{\text{poiss}} = \frac{1}{N_{\text{int}}} \sum_{\mathbf{x}_i \in \Omega} |e^{\pm}_{\text{poiss}}(\mathbf{x}_i)| \tag{11}$$

and analogously for $e^{\max}_{\text{lap}}$ and $e^{\text{avg}}_{\text{lap}}$. $N_{\text{int}}$ is the number of discretisation points inside the domain $\Omega$.

In the next section we will calculate the approximation error for various stencil sizes $n$ and further investigate its (non-trivial) behaviour.

## 3. Core observations

### 3.1. Error oscillations

In Fig. 2 we see that both $e^{\max}_{\text{poiss}}(n)$ and $e^{\text{avg}}_{\text{poiss}}(n)$ oscillate with several local minima (at stencil sizes $n = 28, 46$) and maxima (at stencil sizes $n = 17, 34$). Additionally, these oscillations are not erratic, but instead seem to somewhat resemble a smooth function. This is even more evident in Laplacian approximation error $e^{\text{avg}}_{\text{lap}}(n)$, which is also plotted and we can observe that it has local minima and maxima at same stencil sizes. Such regularity implies that the appearance of these minima in the error is not merely a coincidence, but a consequence of a certain mechanism that could be explained. Further understanding of this mechanism would be beneficial, as it could potentially allow us to predict the location of local minima a priori. Considering that the difference between the neighbouring local maxima and minima can be over an order of magnitude this could greatly increase the accuracy of the method without having to increase the order of the augmentation or the discretisation density. Note that the behaviour of $e^{\max}_{\text{lap}}(n)$ stands out as much more irregular. This implies that in order to explain the observed oscillations, we have to consider the collective behaviour of multiple points, which will be confirmed later on when we consider the error's spatial dependence.

### 3.2. Pointwise behaviour

Fig. 3 provides some more insight into the mechanism behind the oscillating error. Here we have plotted the spatial dependence of the signed error $e^{\pm}_{\text{poiss}}$ for those stencils that correspond to the marked local extrema. We can observe that in the maxima, the error has the same sign throughout the whole domain. On the other hand, near the values of $n$ that correspond to the local minima there are parts of the domain that have differing error signs. Concretely, the sign of $e^{\pm}_{\text{poiss}}$ is negative for stencil sizes between 17 and 27 inclusive. In the minima at $n = 28$ both error signs are present, while for bigger stencil sizes (between 29 and 45 inclusive) the error again has constant sign only this time positive.

This connection between the sign of $e^{\pm}_{\text{poiss}}$ and the minima in $e^{\max}_{\text{poiss}}(n)$ motivates us to define a new quantity:

$$\delta N^{\pm}_{\text{poiss}} = \frac{1}{N_{\text{int}}} \Big( |\{\mathbf{x}_i \in \Omega : e^{\pm}_{\text{poiss}}(\mathbf{x}_i) > 0\}| - |\{\mathbf{x}_i \in \Omega : e^{\pm}_{\text{poiss}}(\mathbf{x}_i) < 0\}| \Big) \tag{12}$$

and analogously for $\delta N^{\pm}_{\text{lap}}$. Simply put, the quantity $\delta N^{\pm}_{\text{poiss}}$ is proportional to the difference between the number of nodes with positively and negatively signed error. Assigning values of $\pm 1$ to positive/negative errors respectively, this quantity can be roughly interpreted as the average sign of the error. It should hold that $|\delta N^{\pm}_{\text{poiss}}|$ is approximately equal to 1 near the maxima and lowers in magnitude as we approach the minima. Fig. 4 confirms this intuition: $\delta N^{\pm}_{\text{poiss}}(n)$ changes its values between $\pm 1$ very abruptly only near the $n$ that correspond to the minima of $e^{\max}_{\text{poiss}}(n)$. A similar conclusion can be made for $\delta N^{\pm}_{\text{lap}}$, which acts as a sort of “smoothed out” version of $\delta N^{\pm}_{\text{poiss}}(n)$.

At a first glance, $\delta N^{\pm}_{\text{lap}}$ seems like a good candidate for an error indicator — it has a well-behaved dependence on $n$, approaches $\pm 1$ as we get closer to the error maxima and has a root near the error minima. The major downside that completely eliminates its applicability in the current state is the fact that we need access to the analytical solution to be able to compute it.

[1] https://gitlab.com/e62Lab/public/2023_cp_iccs_stencil_size_effect

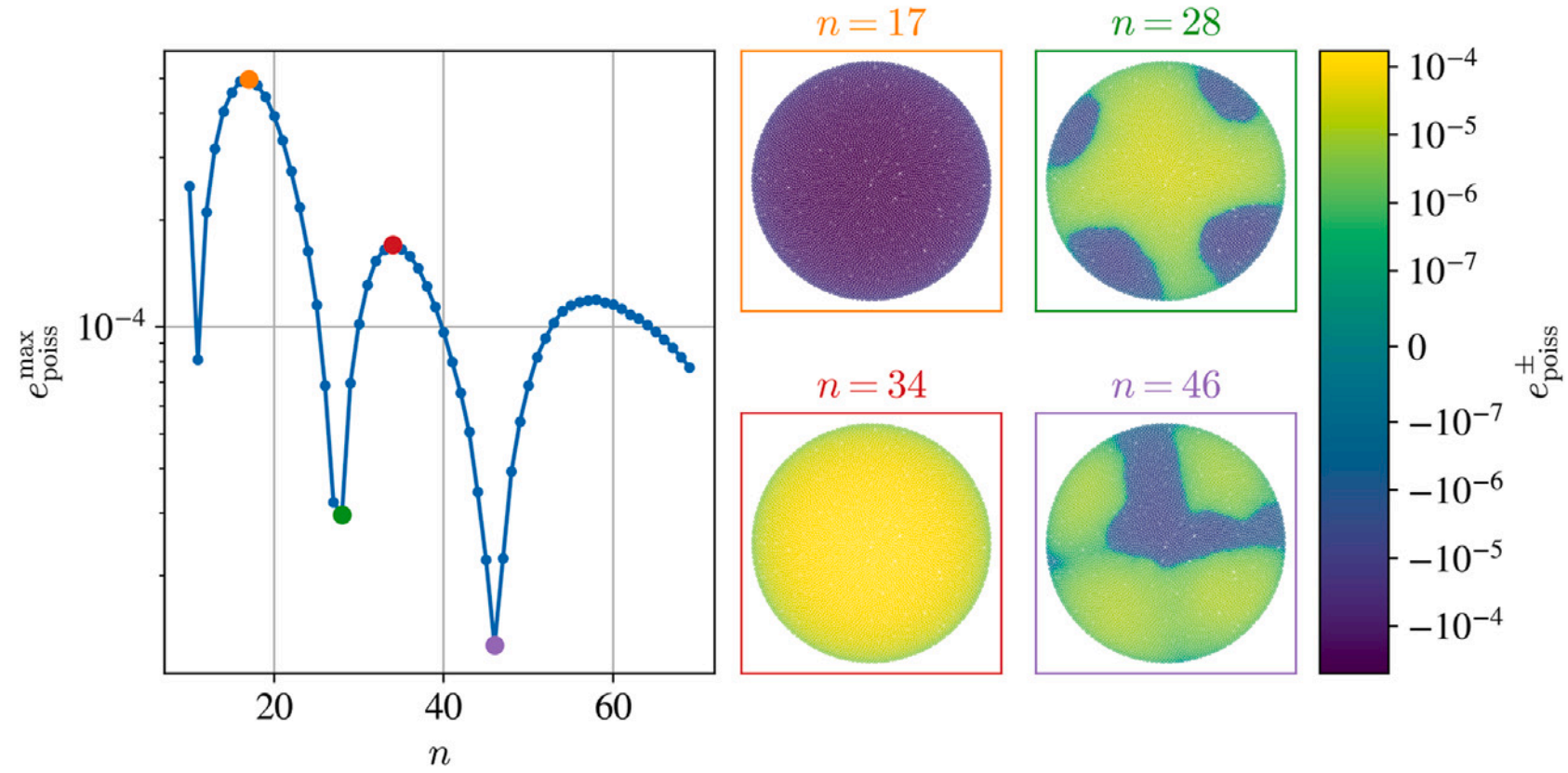

**Fig. 3.** Spatial dependence of $e^{\pm}_{\text{poiss}}$ in some local extrema. The colour scale is the same for all drawn plots.

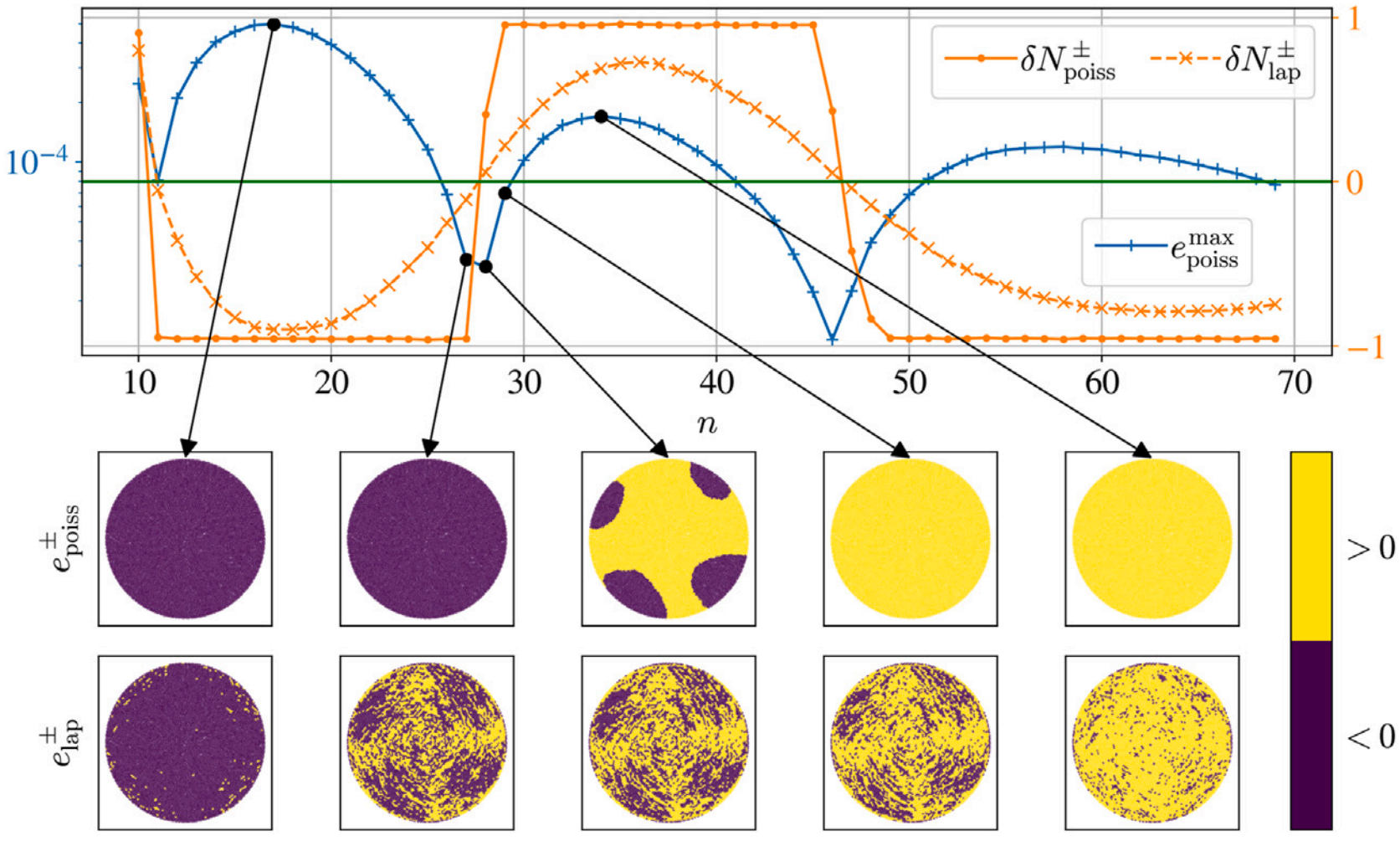

**Fig. 4.** The quantities $\delta N^{\pm}(n)$ along with the spatial profiles of the signs of $e^{\pm}$ for some chosen stencil sizes. For convenience, $\delta N^{\pm} = 0$ is marked with a green line.

## 4. Analysis

Firstly, we would like to eliminate the possibility of some common numerical issues being a cause for the oscillating error. An immediate idea is that the choice of a sparse solver employed at the end of the solution procedure is responsible for the observed behaviour. We have eliminated this possibility by repeating the analysis with both the SparseLU and BiCGSTAB solvers, where no difference has been observed. Likewise, we have checked the condition numbers of the matrices involved in the approximation and have not seen any unexpected behaviour. For further analysis, we will change some aspect of our solution procedure and verify that the observed error oscillations remain.

### 4.1. Changes in the solution procedure

#### 4.1.1. Discretisation refinement

The next idea we explore is the possibility of the discretisation being too coarse. Fig. 5 shows that under discretisation refinement $e^{\max}_{\text{poiss}}(n)$ maintains roughly the same shape and is just shifted vertically towards a lower error. The latter shift is expected, as we are dealing with a convergent method, for which the solution error behaves as $e \propto h^p$ as $h \to 0$, where $p$ is the order of the method, which for the Poisson equation in our setup should equal $p = 2$ [24]. We also show $e^{\max}_{\text{poiss}}(h)$ in a log–log scale for the stencil sizes between $n = 34$ (local maximum) and $n = 46$ (local minimum). It can be seen that the slopes and therefore the orders $p$ do not deviate from the predicted $p = 2$ behaviour as we increase the stencil size and that the observed oscillations mainly come from the proportionality constant in $e \propto h^p$. The stencil dependence of the error proportionality constant has already been observed in similar methods [24,32].

#### 4.1.2. Node layout

As our next test we consider a possibility of a particularly chosen node layout being responsible for the observed behaviour. On Fig. 6 we have repeated our calculations on several different sets of discretisation points, varied by changing the random seed in the DIVG algorithm. Notice that the maximum error oscillations remain almost unchanged. We can observe similar behaviour on discretisations of different character, namely nodes generated by a uniform discretisation in polar coordinates (discretising concentric circles of increasing radii) and by Halton sequences [38]. The case of Halton nodes exhibits considerably less stable behaviour, which is expected as, unlike the other discretisation sets employed, Halton nodes do not have any minimal spacing guarantee. This can be measured by the node set ratio $\gamma$ [39], which is also displayed on the same figure. It is defined as $\gamma = \rho/\delta$, where $\rho = \max_{\mathbf{x}\in\Omega} \min_i \|\mathbf{x}_i - \mathbf{x}\|$ are the fill distance (maximum empty sphere diameter) and $\delta = \min_{i\neq j} \|\mathbf{x}_i - \mathbf{x}_j\|$ the separation distance of a given discretisation set.

#### 4.1.3. Boundary stencils

Next we check if boundary stencils are responsible for the observed behaviour as it is known that they can be problematic due to their one-sideness [30]. On Fig. 7 we have split our domain into two regions

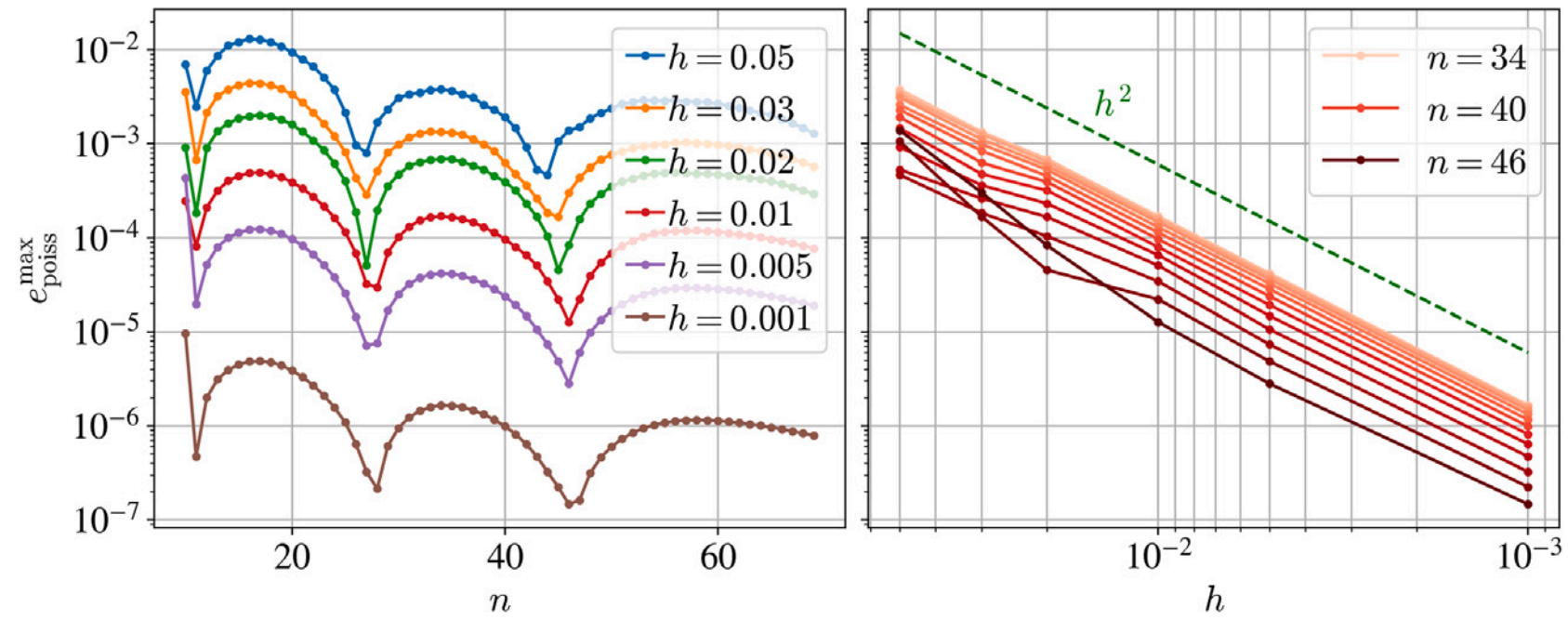

**Fig. 5.** Behaviour of the approximation errors under a refinement of the discretisation.

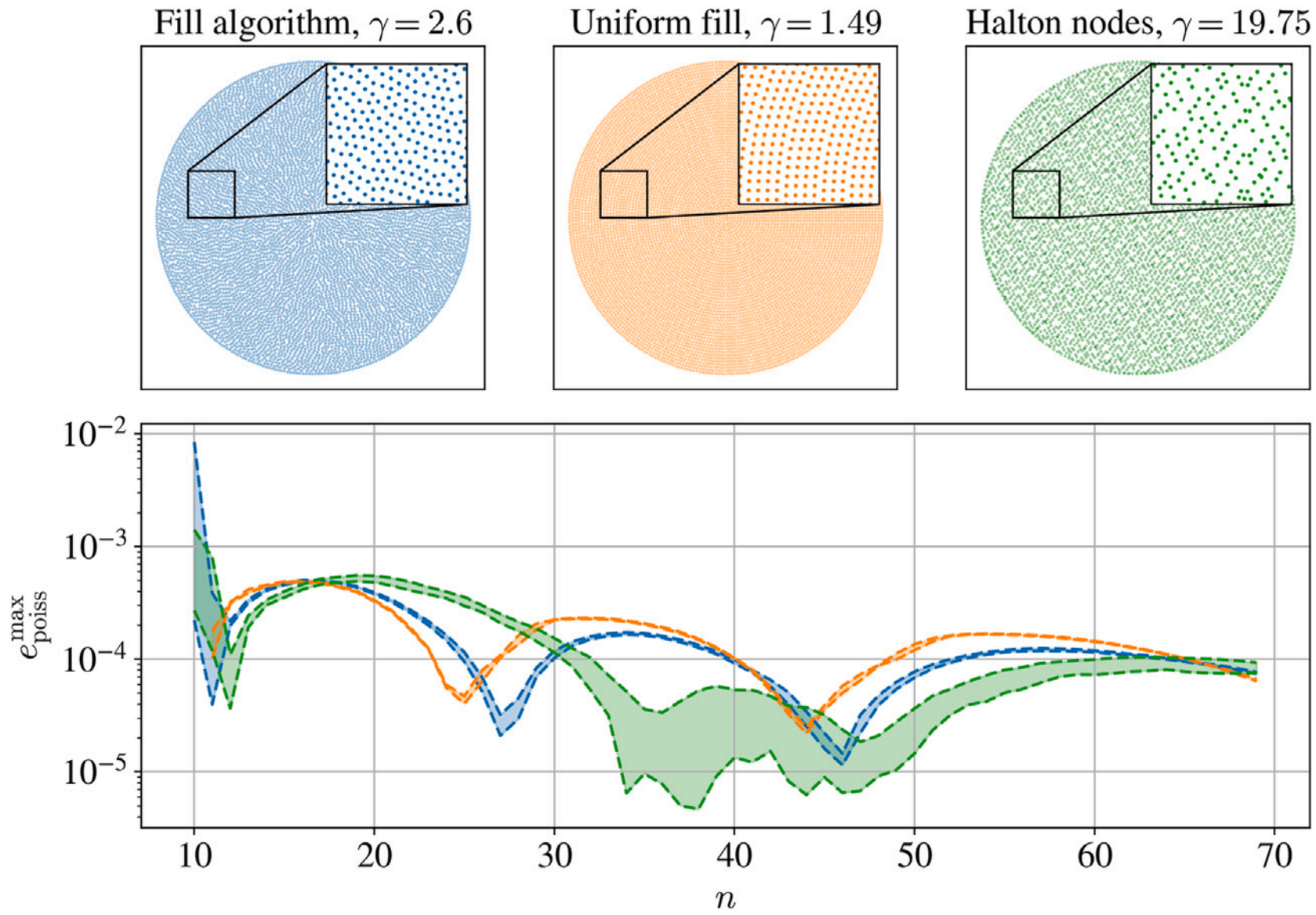

**Fig. 6.** Effect of different discretisation methods on the error oscillations. Example node sets along with their node set ratio $\gamma$ are displayed in the first row. Second row displays the $e^{\max}_{\text{poiss}}$ behaviours: with each of the three listed methods we have generated 10 different node sets and repeated our analysis. The computed values all lie within the corresponding shaded regions.

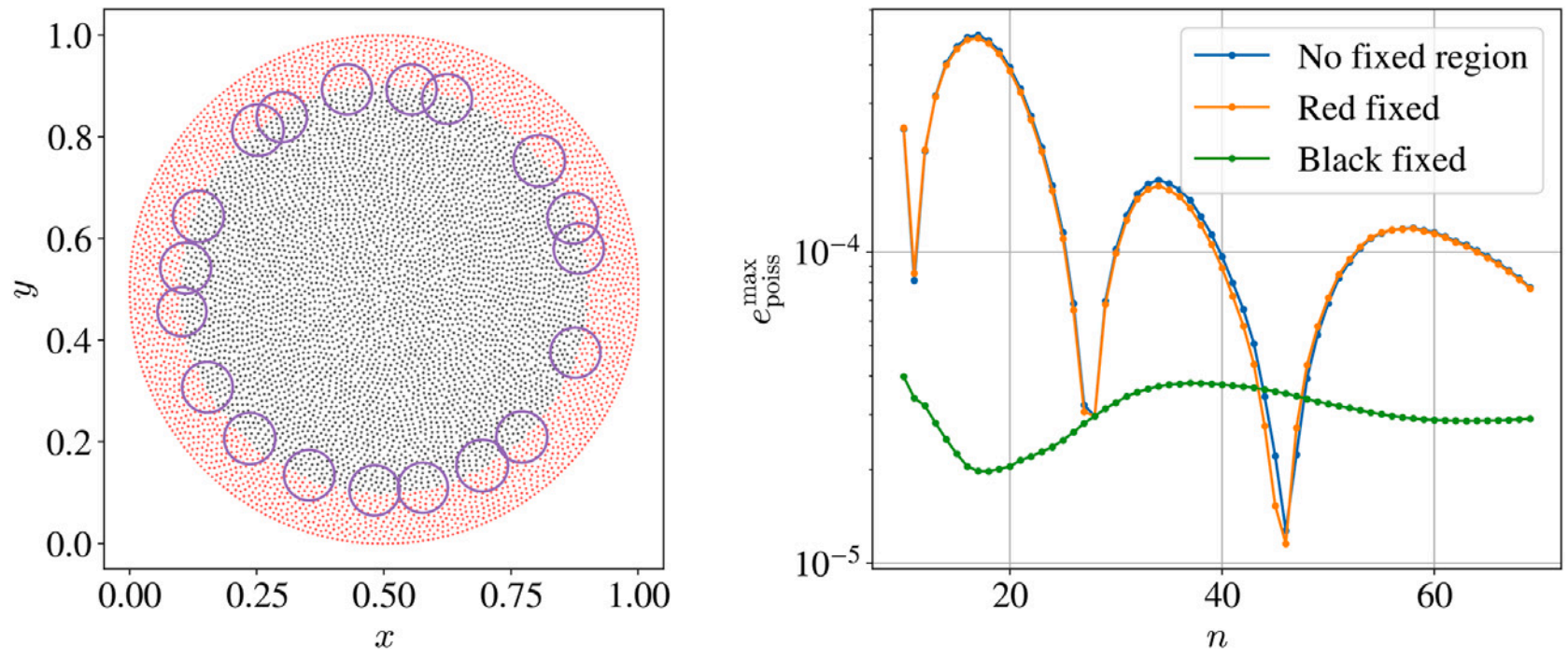

**Fig. 7.** The separation of the domain into two regions is seen on the left, where the purple circles show the radii of the biggest stencils considered ($n = 69$). The right graph shows the error dependence when either of the regions is at a fixed stencil size $n = 28$. The previous result with no fixed stencil size regions is also shown.

— the nodes near the boundary $\{\mathbf{x}_i \in \overline{\Omega} \, : \, \|\mathbf{x}_i - (0.5, 0.5)\| > 0.4\}$ are coloured red, while the nodes far from the boundary $\{\mathbf{x}_i \in \overline{\Omega} \, : \, \|\mathbf{x}_i - (0.5, 0.5)\| \leq 0.4\}$ are black. We can see that the dependence of $e^{\max}_{\text{poiss}}(n)$ marginally changes if we keep the stencil size near the boundary fixed at $n = 28$ (corresponding to one of the previously mentioned minima), while only changing the stencil sizes of the remaining nodes. This shows that the observed phenomenon is not a consequence of the particularly problematic boundary stencils.

#### 4.1.4. Approximation basis

Finally, we check if the oscillations remain also for other choices of parameters associated with our approximation basis. We will consider different types of PHS RBFs: Another radial cubic $\phi(r) = r^{2k+1}$ and also the less commonly used thin plate splines $\phi(r) = r^{2k}\log(r)$. Additionally, for each PHS several different values of the monomial augmentation degree $m$ will be considered, which should not be chosen lower than the PHS parameter $k$ for stability reasons [19]. We should also emphasise that the monomial augmentation degree $m$ affects the minimal possible stencil size so some change in the behaviour is expected. The results can be seen on Fig. 8, where we can see that the stencil size dependence does change as we modify these parameters, but oscillatory behaviour remains.

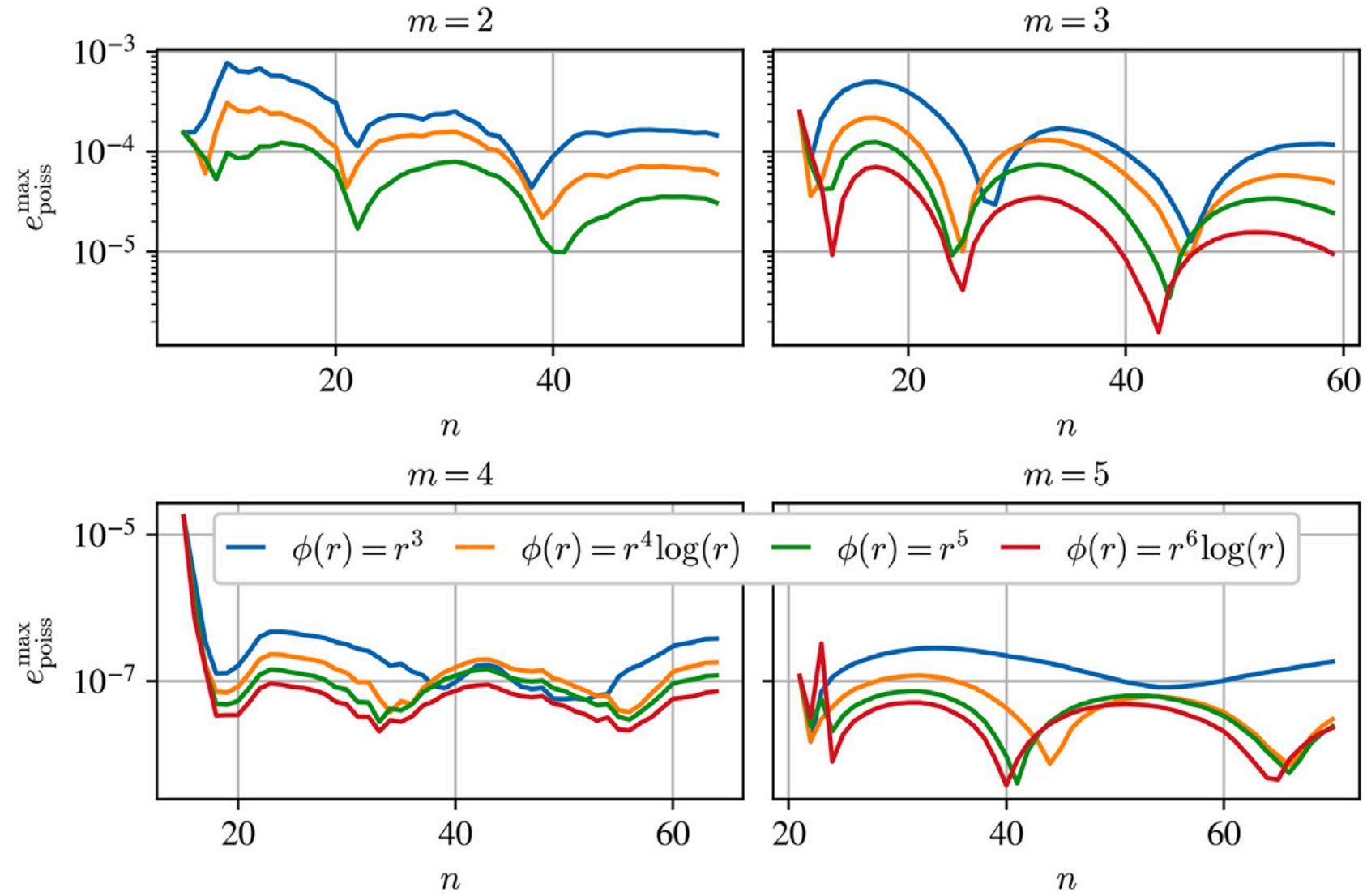

**Fig. 8.** Error oscillations for different choices of method parameters.

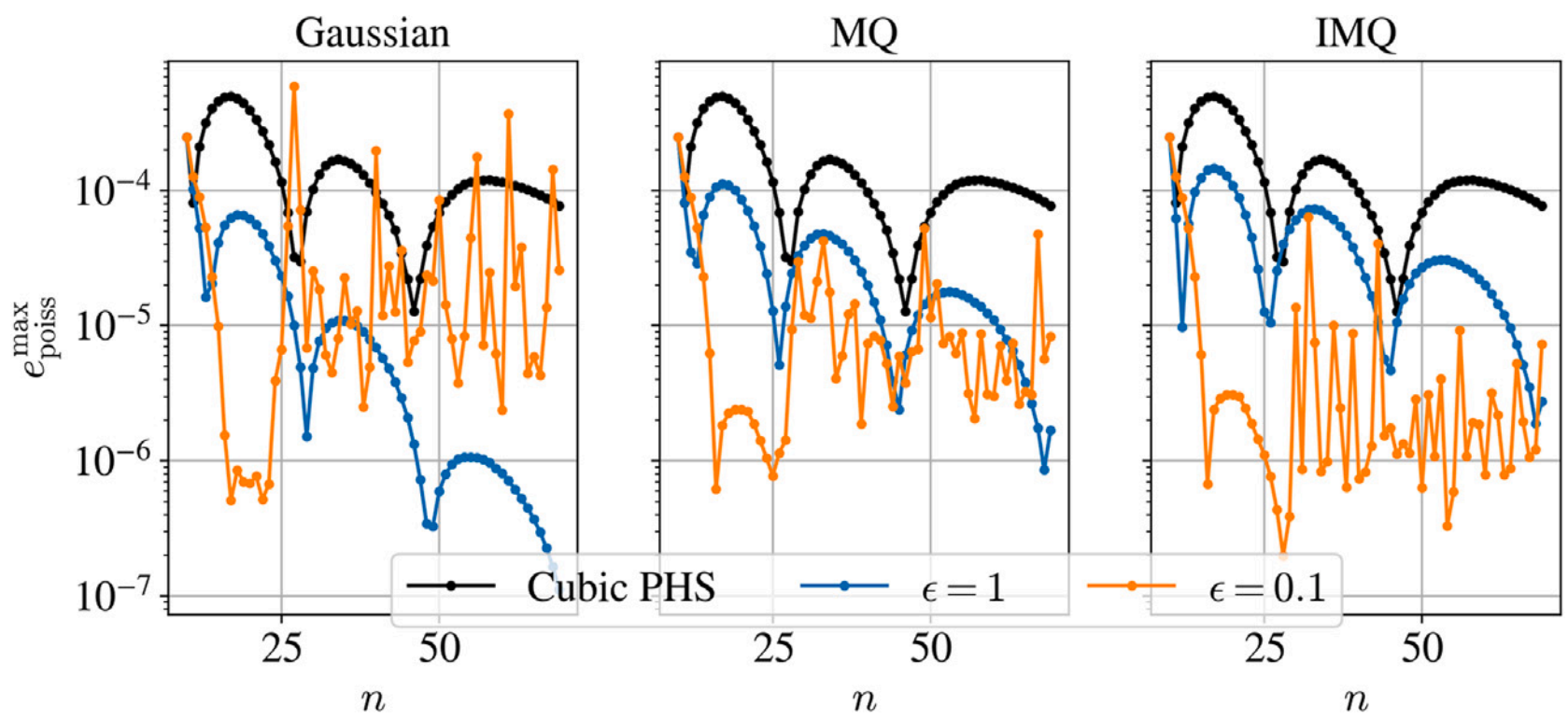

**Fig. 9.** Maximum error dependence on stencil size for few different choices of RBFs. The previous result with $\phi(r) = r^3$ is also shown for reference.

For completeness, let us also briefly consider alternative choices of RBFs besides PHS. We consider the following:

- Gaussian $\phi(r) = e^{-(\epsilon r)^2}$,
- Multiquadric (MQ) $\phi(r) = \sqrt{1 + (\epsilon r)^2}$,
- Inverse Multiquadric (IMQ) $\phi(r) = (1 + (\epsilon r)^2)^{-0.5}$.

We repeat our analysis and study the stencil size dependence $e_{\text{max}}^{\text{poiss}}(n)$. The results can be seen on Fig. 9. We can observe that we do get similar oscillatory behaviour also for these other RBFs when $\epsilon = 1$, but the behaviour is very different when $\epsilon = 0.1$. The fact that $\epsilon$ greatly affects the accuracy is well known and in general, selecting a suitable $\epsilon$ is a difficult problem [22]. Although in our analyses, we have scaled the local coordinate $r$ with the stencil radius, the optimal shape parameter likely still varies with the stencil, which would greatly complicate further investigations. For this reason we focus purely on PHS for the remainder of the paper, as they have no such parameter[2] — one of the main reasons for considering them in the first place.

[2] One could argue that the PHS exponent is also a parameter, however its effect is much less drastic [40].

### 4.2. Different problem setups

Having verified that the observed oscillations remain as we modify our solution procedure, we now turn our attention to different problem setups. As pointwise behaviour, described by $\delta N^{\pm}_{\text{poiss}}$, could have some value in applications we will now include it in our analyses. Concretely, we would like to investigate the extent to which the observed error oscillations and the associated $\delta N^{\pm}_{\text{poiss}}$ behaviour remain if we modify our model problem in some way.

#### 4.2.1. Dimensionality

We start by changing the dimensionality of the domain. Namely, in addition to a 2D disc, we also consider an interval $\Omega = (0, 1)$ and an open ball in 3D, centred around a point $(0.5, 0.5, 0.5)$ with a radius of 0.5, again with the Dirichlet boundary conditions calculated from our chosen analytical solution, which is $u(x) = \sin(\pi x)$ and $u(x, y, z) = \sin(\pi x)\sin(\pi y)\sin(\pi z)$ for the new 1D and 3D cases respectively. The results are displayed on Fig. 10 and we can see that both the error oscillations and the associated pointwise behaviour ($\delta N^{\pm}_{\text{poiss}}$) remain, where the error minima are dimension-dependent and further apart from each other with increasing dimensionality. Furthermore, we can see that the oscillations eventually dampen for sufficiently high stencil sizes. There are also some deviations from our previous pointwise error observations in 2D – as we pass the minimum just below $n = 100$, a notable change can be seen in $\delta N^{\pm}_{\text{poiss}}$ but of different character than previously observed – the sign does not change. Additionally in the 3D case, the value of $|\delta N_{\text{poiss}}^{\pm}|$ is not 1 outside the minima but is instead slightly lower. In spite of these differences, qualitatively the observed behaviour remains the same also in those newly analysed cases, especially at lower and therefore (due to the high computational cost of increasing $n$) the most relevant stencil sizes.

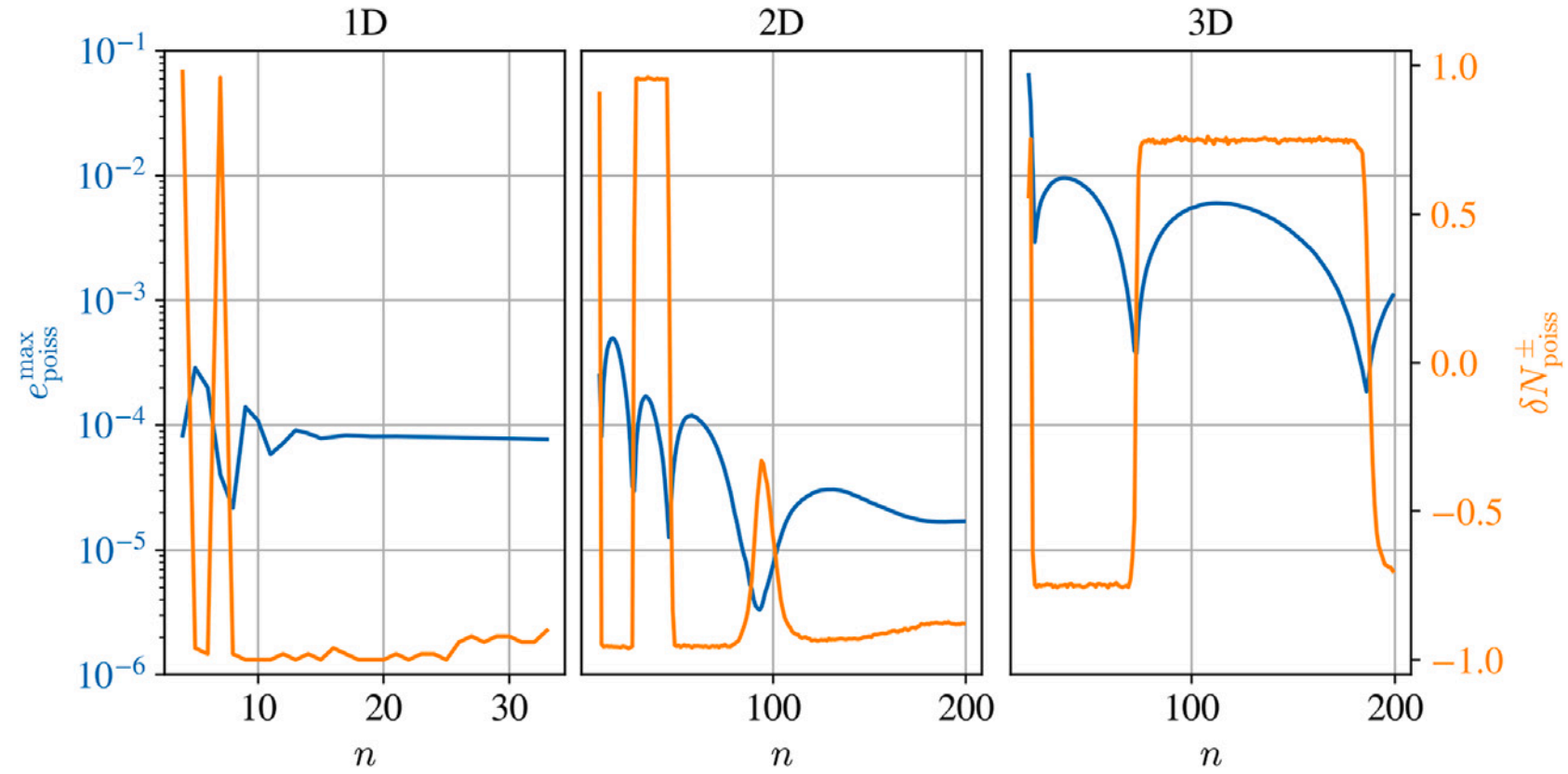

**Fig. 10.** Error oscillations and the $\delta N_{\text{poiss}}^{\pm}$ behaviour plotted for different dimensionalities. The 3D case was computed with a spacing of $h = 0.04$ due to its increased computational complexity.

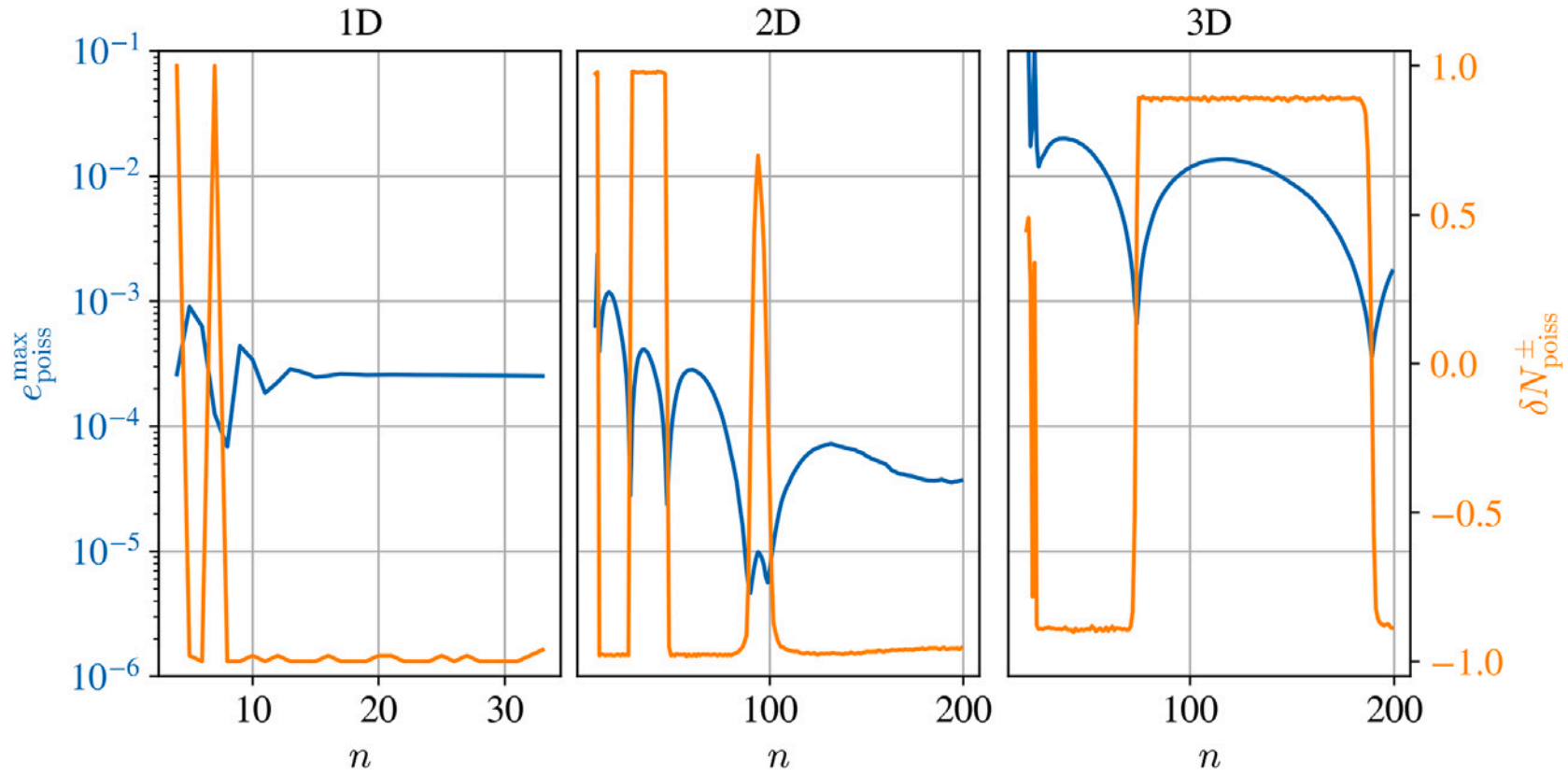

**Fig. 11.** Error oscillations and the $\delta N_{\text{poiss}}^{\pm}$ behaviour plotted for different dimensionalities and with mixed boundary conditions. The 3D case was computed with a spacing of $h = 0.04$ due to its increased computational complexity.

#### 4.2.2. Boundary conditions

We repeat the analysis in different dimensions also with mixed boundary conditions, keeping Dirichlet boundary conditions for boundary nodes satisfying $x > 0.5$ and imposing Neumann boundary conditions for the remaining nodes. In order to increase stability, we have taken a standard approach and added a ghost node for each Neumann boundary node, positioned at a distance $h$ in the outward-facing normal direction [30]. To maintain a square system additional equations must be added — we require the discretised PDE to hold not only at interior points, but also at each Neumann boundary node. The results are seen on Fig. 11. The changes compared to the previously considered case with only Dirichlet boundary conditions are minimal, a notable difference appearing only in the 2D case, where two minima at very close proximity appear just below $n = 100$. This particular minimum already exhibited unusual behaviour previously, on Fig. 10.

#### 4.2.3. Domain shapes

Next, we solve the Poisson problem on different domains, while maintaining the Dirichlet boundary conditions and the same right hand side $f(x, y)$. We have tested the following domain shapes:

- A non-Lipschitz domain — a nephroid, parametrically given by $x(t) = 0.75\cos(t) - 0.5\cos^3(t)$ and $y(t) = 0.5\sin^3(t)$, $t \in [0, 2\pi)$ and then rotated by $\pi/4$ and translated by $(0.5, 0.5)$.
- A triangular domain — a polygon with vertices $(0, 0)$, $(1, 0)$ and $(0, 1)$.
- "Pac-Man" shaped domain — the set-theoretic difference of a disc of radius 0.5 centred at $(0.5, 0.5)$ and a square with corners $(0.5, 0)$ and $(1, 0.5)$. The resulting shape is rotated by $\pi/4$.
- A domain that is not simply connected — an annulus, centred around $(0.5, 0.5)$ with inner radius 0.1 and outer radius 0.5.
- A more complicated domain, given in polar coordinates as $r(\phi) = 0.25\ |\cos(1.5\phi)|^{\sin(3\phi)}$ and also translated by $(0.5, 0.5)$.

Fig. 12 shows that again our previous observations hold with a slight difference in the pointwise error behaviour — $\delta N_{\text{poiss}}^{\pm}$ is again not near zero only in the error minima, but also in their neighbourhood. Additionally on some domains $|\delta N_{\text{poiss}}^{\pm}|$ is not exactly equal to 1 far from the minima as before, but is instead slightly lower. It would appear that the observed phenomena are robust under a change of the domain since, again, while minor differences are present, qualitatively the behaviour is still the same.

#### 4.2.4. Differential operators

To further test the extent of observed oscillations, let us look at more PDEs of type $\mathcal{L}u(x, y) = f(x, y)$, still with the same analytical solution and Dirichlet boundary conditions, but different elliptic partial differential operators $\mathcal{L}$:

$$\mathcal{L}_1 = \nabla^2 + \partial_x\partial_y$$
$$\mathcal{L}_2 = \nabla^2 + \partial_x + \partial_y$$
$$\mathcal{L}_3 = x\partial_x^2 + y^2\partial_y^2$$

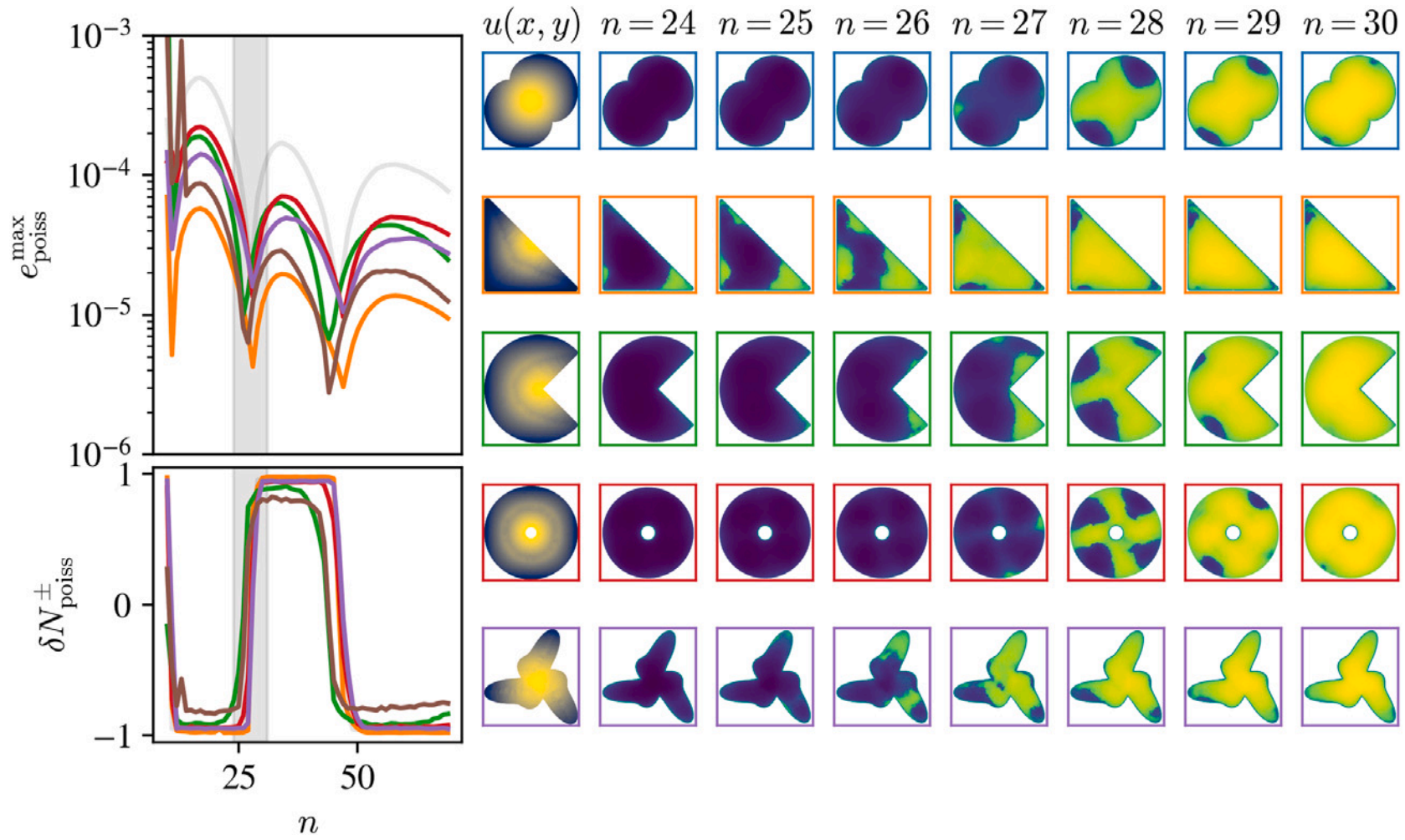

**Fig. 12.** The plots show the error oscillations and its pointwise behaviour for different domains. The shaded region corresponds to the interval $n \in [24, 30]$, for which the pointwise contours are displayed to the right. The contour plot borders' colour associates them to the corresponding plots on the left. For reference, the results of our initially considered setup are plotted with a black, transparent line.

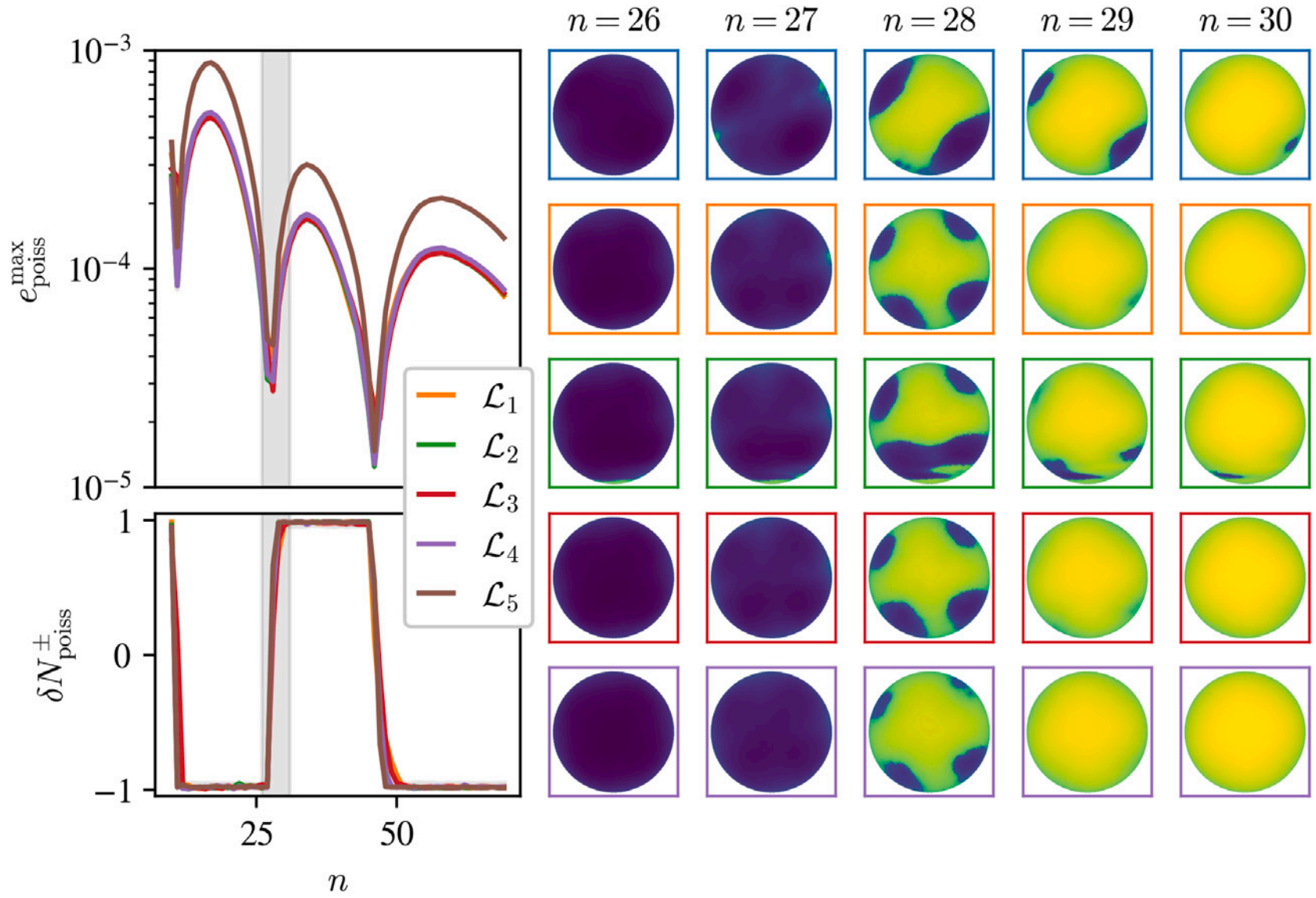

**Fig. 13.** The plots on the left show the $e_{\text{poiss}}^{\max}(n)$ and $\delta N_{\text{poiss}}^{\pm}(n)$ dependence respectively for several different elliptic equations. The shaded region corresponds to the interval $n \in [26, 30]$, for which the pointwise contours are displayed to the right. The contour plot border's colour associates them to the corresponding plots on the left. For reference, the results of our initially considered setup are plotted with a black, transparent line.

$$\mathcal{L}_4 = \nabla^2 + 1$$

$$\mathcal{L}_5 = \nabla^2 + 10$$

The results are displayed on Fig. 13, where we can see that the behaviour again remains very similar with the minima appearing at exactly the same stencil sizes for different differential operators.

#### 4.2.5. Analytical solutions

As a last test of the section, we return to $\mathcal{L} = \nabla^2$ and instead consider other choices of the analytical solution $u(x, y)$. Consequently this changes the right hand side $f(x, y)$ of the Poisson problem as well as its boundary conditions. A variety of different functions were tested, listed in Table 1. The results are shown on Fig. 14. We see that oscillatory behaviour in the error is observed in all of the listed functions, however, unlike in most of our previous tests, locations of the error minima can substantially change and are function-dependent. Additionally, in the case of $u_6$ (the well known Franke function, commonly used for testing method robustness) minima are also harder to see, likely due to the increased complexity of this particular function. Most importantly, we notice that the previously observed pointwise error behaviour is not general, since it differs for the cases of $u_4, u_5$ and $u_6$. Nevertheless, a visual connection can be still made in all the cases: $u_5$ has $\delta N_{\text{poiss}}^{\pm}(n)$ approximately equal to zero, with bumps at the stencil sizes corresponding to minima locations. The remaining two cases actually seem to exhibit exactly the opposite of our initial observations — in minima the pointwise errors are of mostly one sign while in maxima both signs are approximately equally represented.

**Table 1**
The different choices of analytical solution $u(x, y)$ considered.

| Label | $u(x, y)$ |
|---|---|
| $u_1$ | $x^4 y^5$ |
| $u_2$ | $1 + \sin(4x) + \cos(3x) + \sin(2y)$ |
| $u_3$ | $\exp(x^2)$ |
| $u_4$ | $\operatorname{arsinh}(x + 2y)$ |
| $u_5$ | $\cos(\pi x)\cos(\pi y)$ |
| $u_6$ | $\operatorname{franke}(x, y)$ [41] |

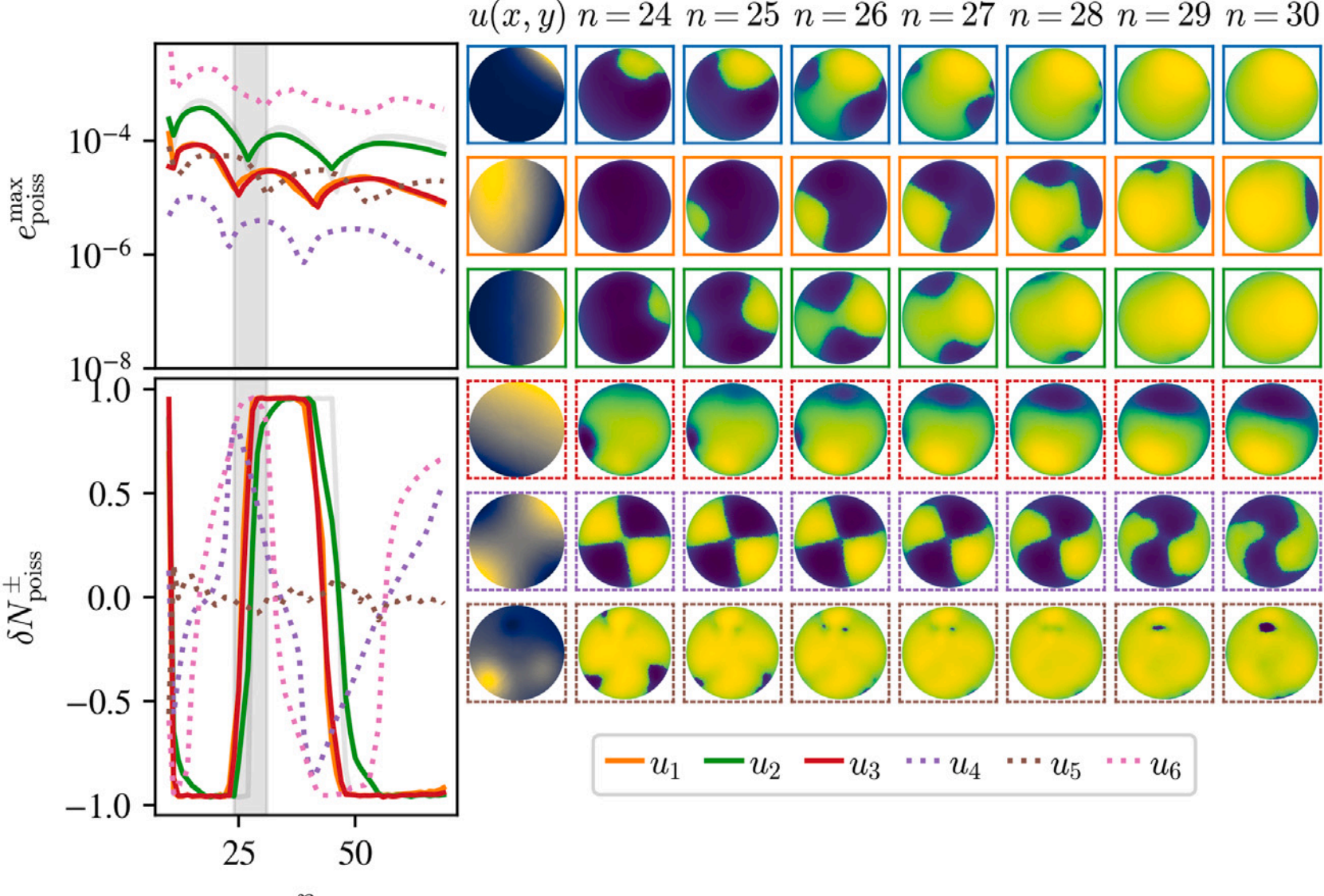

**Fig. 14.** The plots on the left show the $e^{\max}_{\text{poiss}}(n)$ and $\delta N^{\pm}_{\text{poiss}}(n)$ dependence respectively for several different choices of solution functions. The shaded region corresponds to the interval $n \in [24, 30]$, for which the pointwise contours are displayed to the right. The contour plot border's colour associates them to the corresponding plots on the left. Graphs plotted with a dotted line correspond to functions of which pointwise error behaviour differs from previously observed. For reference, the results of our initially considered setup are plotted with a black, transparent line.

It could therefore still be possible to deduce the location of error minima purely from $\delta N^{\pm}_{\text{poiss}}(n)$. However, in order to do so we need to understand which properties of $u(x, y)$ (or, preferably, easily obtainable properties of the Poisson problem) determine the three observed pointwise error behaviours. We have thus far not explored this further.

## 5. A practical example

We conclude the paper with a practical example, hinting at a potential application of the presented results. We will consider the problem of determining the steady state temperature profile of a heatsink, similarly as in [42]. As a 3D domain with irregular boundary, this is a typical case where meshless methods can be employed. The domain description has been obtained from [43]. Converting our physical problem to the context of PDEs, determining the steady state temperature profile amounts to solving the Laplace equation:

$$\nabla^2 T(x, y) = 0. \tag{13}$$

At the bottom (minimal $y$), the heatsink will be in contact with a body of a temperature $T_{\text{hot}}$, giving a Dirichlet boundary condition. On the rest of the boundary we can combine Fourier's law of conduction and Newton's law of cooling to obtain the Robin boundary condition $T + \alpha \frac{\partial T}{\partial n} = T_{\text{out}}$, where $\alpha = \lambda / h$, with $\lambda$ being the thermal conductivity of the heatsink and $h$ the heat transfer coefficient. $T_{\text{out}}$ is the ambient temperature.

To fully specify the problem setup, we will set the temperatures equal to $T_{\text{out}} = 20\,^{\circ}\text{C}$, $T_{\text{hot}} = 80\,^{\circ}\text{C}$. Heatsinks are commonly made from aluminium with $\lambda \approx 209\,\text{W}\,\text{m}^{-1}\,\text{K}^{-1}$. Outside the heatsink, we have air moving at a moderate speed, which amounts to $h \approx 100\,\text{W}\,\text{m}^{-2}\,\text{K}^{-1}$.

Due to the higher complexity of this problem, the parameters were modified from the ones most commonly used until now. We now work with $h = 3.5 \times 10^{-4}$ m and $m = 2$. Example numerical solution can be seen on the left plot of Fig. 15.

As we do not have an analytical solution for this case, we have resorted to estimating the behaviour of solution error with the recently introduced IMEX error indicator [10,44]. IMEX uses an auxiliary higher order operator to obtain an estimation of error behaviour and works well provided the solution is smooth enough (as it appears to be also in our case). For more information on IMEX we refer an interested reader to the paper [10].

Right plot of Fig. 15 shows the IMEX (calculated pointwise and averaged over the domain) dependence on stencil size. Once again we observe oscillations that are in agreement with our previous observations — the minima are quite far apart (we are in a 3D case), while the difference between the extrema can be substantial. As already mentioned, we currently do not have all the building blocks to implement the presented application in practice, as we do not have a systematic way of obtaining the optimal stencil size without access to some accurate reference solution. This presented example is only meant to provide additional motivation for further research in this direction. Additionally, it should be noted that we have two possible methods of increasing the accuracy of our solution procedure — lowering the discretisation distance $h$ or suitably modifying the stencil size $n$. There is likely some trade-off here where for some cases each of the options turn out to be superior, but we have not delved further into the subject.

## 6. Conclusions

Our study started with a simple observation — when solving a Poisson problem with PHS RBF-FD, the approximation accuracy depends on the stencil size $n$ in a non-trivial manner. In particular, there exist certain stencil sizes where the method is especially accurate. A priori knowledge of these stencil sizes could decrease the solution error without any additional effort and is therefore strongly desirable. We have made a small step towards understanding this phenomenon by looking at the spatial dependence of the signed solution error — in the stencil sizes corresponding to the local error minima, the signed solution error was not strictly positive or negative. This was unlike the generic stencil sizes, where the error had the same sign throughout the domain. Motivated by this observation, we have introduced a quantity that is roughly the average sign of the pointwise errors and appeared to have a root only near the stencils corresponding to the local error minima.

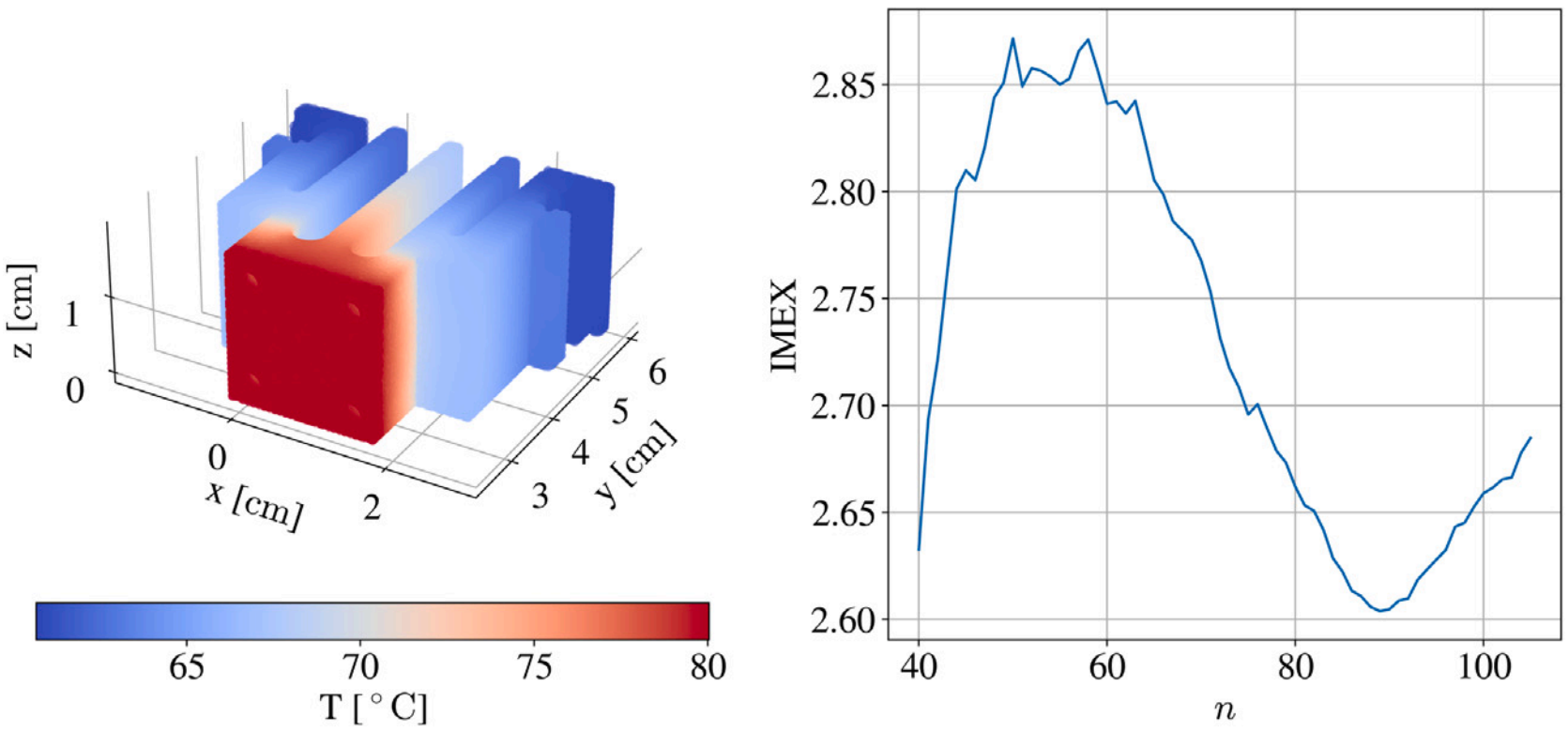

**Fig. 15.** Left plot shows the numerical solution to the physical problem, calculated with the stencil size of $n = 90$. Right plot shows the IMEX dependence on stencil size.

We have eliminated some common numerical issues as the cause of the observed oscillations and verified that the latter remain also when changing some aspect of our numerical solution procedure. This was followed by extending our analyses to different problem setups, where we have noticed that the observed behaviour is robust under the change of the problem domain or its dimensionality and is not limited to the simplest Dirichlet boundary conditions. Furthermore the behaviour remains almost unchanged also when considering some other PDEs of elliptic type. However, we have observed a greater effect of the chosen analytical solution (which in turn affects the right hand side of the Poisson problem as well as the boundary conditions). While oscillatory behaviour of the aggregated error remains, its pointwise behaviour can differ. Nevertheless, we can still hope to be able to deduce the location of minima from the pointwise behaviour, since we have seen that a connection between the maximum error and its spatial properties can still be made, although this relationship is more complicated than initially described and warrants further investigation.

As demonstrated by the heatsink example, our observations are not limited to simple, unrealistic problems, but could extend to real-life examples as well. The research presented is a step towards defining a more practically useful indicator, which would reveal the most accurate stencil sizes even without having access to the analytical solution. Such an indicator would provide us with a simple way of improving our method's accuracy — by modifying the stencil size. This represents an alternative to increasing the method's order or refining the discretisation and is part of our ongoing research. Additional future work includes more rigorous theoretical explanations for the observations presented, especially on characterising the different pointwise error behaviours. The effect of different method parameters was merely touched upon in this paper and should also be more thoroughly explored. Further experimental investigations should also be made, particularly to what extent our observations carry over to even more problem setups — differential equations that are not of the elliptic type and a sweep over an even wider choice of analytical solutions, for instance.

## CRediT authorship contribution statement

**Andrej Kolar-Požun:** Writing – review & editing, Writing – original draft, Visualization, Software, Methodology, Investigation, Formal analysis, Conceptualization. **Mitja Jančič:** Writing – review & editing, Investigation, Conceptualization. **Miha Rot:** Writing – review & editing, Conceptualization. **Gregor Kosec:** Project administration, Methodology, Funding acquisition, Conceptualization, Supervision, Writing – review & editing.

## Declaration of competing interest

The authors declare that they have no known competing financial interests or personal relationships that could have appeared to influence the work reported in this paper.

## Acknowledgements

The authors acknowledge the financial support from the Slovenian Research and Innovation Agency (ARIS) in the framework of the research core funding No. P2-0095, the Young Researcher programs PR-10468 and PR-12347, and research projects No. J2-3048 and No. N2-0275.

Funded by National Science Centre, Poland under the OPUS call in the Weave programme 2021/43/I/ST3/00228. This research was funded in whole or in part by National Science Centre (2021/43/I/ST3/00228). For the purpose of Open Access, the author has applied a CC-BY public copyright licence to any Author Accepted Manuscript (AAM) version arising from this submission.

## Appendix. Tabulated error values

For the first case considered, where we have initially demonstrated the phenomenon of oscillating error (Fig. 2), the exact error values are listed in Table A.2. The error values for other cases are omitted and we again direct an interested reader to our public gitlab repository.[3]

**Table A.2**
Error values for the case at the beginning of Section 3, seen on Fig. 2.

| $n$ | $e^{poiss}_{max}$ | $e^{poiss}_{avg}$ | $e^{lap}_{max}$ | $e^{lap}_{avg}$ |
|---|---|---|---|---|
| 10 | 0.000248 | 7.2e−05 | 3.44 | 0.00397 |
| 11 | 8.1e−05 | 2.46e−05 | 0.0473 | 0.00279 |
| 12 | 0.00021 | 8.12e−05 | 0.0397 | 0.00329 |
| 13 | 0.000317 | 0.000129 | 0.02 | 0.00399 |
| 14 | 0.000404 | 0.000169 | 0.0208 | 0.00469 |
| 15 | 0.000457 | 0.000193 | 0.0187 | 0.00516 |
| 16 | 0.000491 | 0.000209 | 0.0203 | 0.00546 |
| 17 | 0.000498 | 0.000213 | 0.02 | 0.00553 |
| 18 | 0.000479 | 0.000206 | 0.0206 | 0.00533 |
| 19 | 0.000443 | 0.000191 | 0.0189 | 0.00498 |
| 20 | 0.000392 | 0.00017 | 0.0192 | 0.0045 |
| 21 | 0.000335 | 0.000146 | 0.0173 | 0.00398 |
| 22 | 0.000274 | 0.000121 | 0.0165 | 0.00348 |
| 23 | 0.000217 | 9.73e−05 | 0.017 | 0.00306 |
| 24 | 0.000163 | 7.55e−05 | 0.0156 | 0.00273 |
| 25 | 0.000115 | 5.49e−05 | 0.0164 | 0.00245 |
| 26 | 6.84e−05 | 3.35e−05 | 0.0135 | 0.00226 |
| 27 | 3.22e−05 | 1.36e−05 | 0.0115 | 0.00214 |
| 28 | 2.96e−05 | 8.13e−06 | 0.0102 | 0.00208 |
| 29 | 6.96e−05 | 2.44e−05 | 0.0116 | 0.00208 |
| 30 | 0.000102 | 3.95e−05 | 0.011 | 0.00213 |

**Table A.2** (*continued*).

| $n$ | $e_{\max}^{\text{poiss}}$ | $e_{\text{avg}}^{\text{poiss}}$ | $e_{\max}^{\text{lap}}$ | $e_{\text{avg}}^{\text{lap}}$ |
|---|---|---|---|---|
| 31 | 0.000131 | 5.3e−05 | 0.0102 | 0.00219 |
| 32 | 0.000153 | 6.33e−05 | 0.0101 | 0.00223 |
| 33 | 0.000165 | 6.94e−05 | 0.0112 | 0.00227 |
| 34 | 0.00017 | 7.23e−05 | 0.0107 | 0.00226 |
| 35 | 0.000166 | 7.15e−05 | 0.00982 | 0.0022 |
| 36 | 0.000158 | 6.89e−05 | 0.00917 | 0.00211 |
| 37 | 0.000146 | 6.46e−05 | 0.0081 | 0.002 |
| 38 | 0.00013 | 5.82e−05 | 0.00879 | 0.00186 |
| 39 | 0.000114 | 5.19e−05 | 0.00893 | 0.00174 |
| 40 | 9.66e−05 | 4.5e−05 | 0.00781 | 0.00161 |
| 41 | 7.97e−05 | 3.81e−05 | 0.00769 | 0.00149 |
| 42 | 6.53e−05 | 3.17e−05 | 0.00813 | 0.0014 |
| 43 | 5.07e−05 | 2.51e−05 | 0.00806 | 0.00133 |
| 44 | 3.43e−05 | 1.76e−05 | 0.00737 | 0.00126 |
| 45 | 2.21e−05 | 1.02e−05 | 0.00814 | 0.00123 |
| 46 | 1.27e−05 | 3.67e−06 | 0.0147 | 0.00122 |
| 47 | 2.23e−05 | 5.99e−06 | 0.0147 | 0.00122 |
| 48 | 3.92e−05 | 1.23e−05 | 0.014 | 0.00124 |
| 49 | 5.42e−05 | 1.91e−05 | 0.0131 | 0.00127 |
| 50 | 6.84e−05 | 2.55e−05 | 0.013 | 0.00133 |
| 51 | 8.23e−05 | 3.16e−05 | 0.0118 | 0.00137 |
| 52 | 9.29e−05 | 3.68e−05 | 0.0125 | 0.00143 |
| 53 | 0.000103 | 4.15e−05 | 0.0122 | 0.00147 |
| 54 | 0.000111 | 4.52e−05 | 0.0123 | 0.00151 |
| 55 | 0.000115 | 4.76e−05 | 0.012 | 0.00153 |
| 56 | 0.000118 | 4.94e−05 | 0.0112 | 0.00154 |
| 57 | 0.000119 | 5.03e−05 | 0.0128 | 0.00156 |
| 58 | 0.000119 | 5.1e−05 | 0.0139 | 0.00156 |
| 59 | 0.000117 | 5.06e−05 | 0.0142 | 0.00155 |
| 60 | 0.000116 | 5.02e−05 | 0.0138 | 0.00154 |
| 61 | 0.000112 | 4.92e−05 | 0.0134 | 0.00151 |
| 62 | 0.000108 | 4.8e−05 | 0.0132 | 0.00148 |
| 63 | 0.000105 | 4.67e−05 | 0.0181 | 0.00145 |
| 64 | 0.000101 | 4.51e−05 | 0.0182 | 0.00142 |
| 65 | 9.68e−05 | 4.34e−05 | 0.0167 | 0.00138 |
| 66 | 9.22e−05 | 4.15e−05 | 0.0159 | 0.00133 |
| 67 | 8.74e−05 | 3.96e−05 | 0.0167 | 0.00129 |
| 68 | 8.23e−05 | 3.75e−05 | 0.0171 | 0.00124 |
| 69 | 7.69e−05 | 3.52e−05 | 0.017 | 0.00118 |

[3] https://gitlab.com/e62Lab/public/2023_cp_iccs_stencil_size_effect